\documentclass[12pt]{article}
\usepackage{amssymb}

\usepackage{mitpress}

\newtheorem{theorem}{Theorem}

\newtheorem{corollary}[theorem]{Corollary}

\newtheorem{definition}[theorem]{Definition}

\newtheorem{lemma}[theorem]{Lemma}

\newtheorem{proposition}[theorem]{Proposition}

\newenvironment{proof}[1][Proof]{\textbf{#1.} }{\ \rule{0.5em}{0.5em}}
\newdimen\dummy
\dummy=\oddsidemargin
\input{tcilatex}

\begin{document}

\title{Completing Artin's braid group on infinitely many strands}
\author{Paul Fabel \\
Department of Mathematics \& Statistics\\
Mississippi State University}
\maketitle

\begin{abstract}
A generalization of the topological fundamental group is developed in order
to construct a completion of Artin's braid group on infinitely many strands
with respect to the following notion of convergence: $b_{n}\rightarrow id$
iff for each $M>0,$ eventually the first $M$ strands of $b_{n}$ are trivial.
\end{abstract}

\section{Introduction}

Artin's braid group $B_{n}$ is the fundamental group of the collection of
planar sets each of which has $n$ elements. Fixing a basepoint $Z_{n}\subset
R^{2},$ elements of $B_{n}$ can be seen as collections of disjoint arcs in $%
R^{2}\times \lbrack 0,1]$ connecting $Z_{n}\times \{0\}$ to $Z_{n}\times
\{1\}$.

Artin's braid group on infinitely many strands$\,B_{\infty }$ is the direct
limit of $B_{n}$ under monomorphisms $i_{n}:B_{n}\rightarrow B_{n+1}$ that
attach a `trivial strand'. It has been studied from various perspectives by
Broto \cite{brot}, Dehornoy \cite{Deh1} \cite{Deh2}, Dynnikov \cite{DY},
Moran \cite{Moran}, Kim and Rolfsen \cite{Rolfsen}, Vershinin \cite{Versh2}
and others.

The group $B_{\infty }$ becomes a \textit{topological} group under the
requirement that a sequence of braids $b_{n}\rightarrow id\in B_{\infty }$
iff for each $M$ there exists $N$ such that if $n\geq N$ then the first $M$
strands of $b_{n}$ are trivial.

Unfortunately $B_{\infty }$ admits no complete metric compatible with its
topology. The goal of this paper is to enlarge $B_{\infty }$ to a complete
group of braids on infinitely many strands.

Recent work of Biss \cite{Biss} shows generally that the based fundamental
group of a space $\pi _{1}(X,p)$ inherits a canonical topology from $X$ and
becomes a topological group.

In an effort order to complete $B_{\infty },$ we develop a generalization of
the topological fundamental group and subsequently exhibit a completely
metrizable topological group $\overline{B_{\infty }}$ such that $B_{\infty }$
is embedded in $\overline{B_{\infty }}$ as a dense subgroup. By construction
each element of $\overline{B_{\infty }}$ can be seen as a braid with
infinitely many strands, a countable collection of pairwise disjoint arcs in 
$R^{2}\times \lbrack 0,1].$

The bulk of this paper is devoted to proving, for a certain space $X$ ( the
space of injections of the integers into the plane), that the topological
fundamental group of $X$ is isomorphic to $\lim_{\leftarrow }G_{n},$ the
inverse limit of Artin's pure braid groups on $n$ strands. This serves to
complete Artin's \textit{pure} braid group on infinitely many strands, the
subgroup of braids where each strand has identical endpoints. A more general
construction is then employed to handle complications created by permutation
of the strand endpoints.

The construction of $\overline{B_{\infty }}$ in this paper serves as a
counterpart to that of \cite{Fabel} in which $B_{\infty }$ is completed into
a group of tame braids via the mapping class group of a disk with infinitely
many punctures. However the topology of $B_{\infty }$ in \cite{Fabel} is
finer than that of the current paper. In \cite{Fabel} a sequence of braids $%
b_{n}$ converges to $id$ iff the first $M$ strands are eventually trivial 
\textit{and unlinked }from the remaining strands. The latter notion of
convergence in $B_{\infty }$ is compatible with a left order topology of
Dehornoy \cite{Deh1}, and the construction in \cite{Fabel} seems to resolve
the issue of its completion as suggested by Dehornoy in \cite{Deh2}.

By contrast, with the coarser topology of $B_{\infty }$ adopted in the
current paper, some elements of $\overline{B_{\infty }}$ are `wild' in a
strong sense: each representative of a wild element is not the restriction
of an ambient isotopy of the plane.

\section{\label{defs}Preliminaries}

Spaces of continuous functions will have the compact open topology, product
spaces will have the product topology, and quotient spaces will have the
quotient topology.

Given topological spaces $A$ and $Y$ let $C(A,Y)$ denote the continuous
functions from $A$ into $Y.$ Let $I(A,Y)=\{f\in C(A,Y)|$ $f$ is one to one$%
\}.$ Note elements of $I(A,Y)$ are not necessarily embeddings. Let $%
E(A,Y)=\{f\in C(A,Y)|f$ is a surjection$\}.$ Let $H(A,A)$ denote the group
of homeomorphisms from $A$ onto $A.$

A surjective map $q:A\rightarrow Y$ is a \textbf{quotient map} provided $V$
is open in $Y$ if and only if $q^{-1}(V)$ is open in $A.$ The following
propositions (Theorem 11.1 p. 139 and Corollary 5.4 p287) are proved in
Munkres \cite{Munkres}.

\begin{proposition}
\label{qmap}If $q:A\rightarrow Y$ is a quotient map and $f:A\rightarrow B$ a
map such that $f$ is constant on sets of the form $q^{-1}(y)$ then there
exists a map $F:Y\rightarrow B$ such that $f=F(q).$
\end{proposition}

\begin{proposition}
\label{munkcor}If $X$ is locally compact and Hausdorf then $F:X\times
Z\rightarrow Y$ is continuous iff $\stackrel{\symbol{94}}{F}:Z\rightarrow
C(X,Y)$ is continuous where $\stackrel{\symbol{94}}{F}$ is defined by $(%
\stackrel{\symbol{94}}{F}(z))(x)=F(x,z).$
\end{proposition}

If $A_{1},A_{2},...$ is a sequence of sets and $\psi _{n}:A_{n}\rightarrow
A_{n-1}$ is a surjection then the \textbf{inverse limit} of $A_{n}$ under $%
\psi _{n}$ is the denoted $\lim_{\leftarrow }A_{n}$ and consists of all
sequences $(a_{1},a_{2},...)$ such that $a_{n}\in A_{n}$ and $\psi
_{n}(a_{n})=a_{n-1}.$

If the topological space $G$ is a group with multiplication $m$ and
inversion $v$ then $G$ is a \textbf{topological group} provided $m$ and $v$
are continuous.

Let $Z^{+}=\{(n,0)\in R^{2}|n\in \{1,2,3,..\}\}.$

Let $Z_{n}=\{(k,0)\in Z^{+}|1\leq k\leq n\}.$

\begin{lemma}
\label{fib}The map $\pi :H(R^{2},R^{2})\rightarrow I(Z_{n},R^{2})$ defined
via $\pi (h)=h_{|Z_{n}}$ determines a fibre bundle.
\end{lemma}

\begin{proof}
Given $x$ in the interior of a closed round planar disk $D$ there is a
circle's worth of line segments each of which has one endpoint $x$ and the
other endpoint on $\partial D.$ Given $\{x,y\}\subset int(D)$ there is a
canonical homeomorphism $h(D,x,y):D\rightarrow D$ mapping $x$ to $y,$
mapping one segment linearly onto another, and fixing pointwise $\partial D.$
Given $\alpha \in I(Z_{n},R^{2})$ choose $h_{\alpha }$ such that $\pi
(h_{\alpha })=\alpha =h_{\alpha |Z_{n}}.$ Choose disjoint closed disks $%
D_{i} $ centered at $\alpha (i).$ For $\beta $ close to $\alpha $ let $%
h_{\alpha ,\beta }:R^{2}\rightarrow R^{2}$ map $\alpha (i)$ to $\beta (i)$,
leaving invariant $D_{i}$ via the canonical homeomorphism, and fixing
pointwise $R^{2}\backslash \cup D_{i}.$ Let $F=\{h\in H(R^{2},R^{2})|$ $%
h_{|Z_{n}}=id_{|Z_{n}}\}.$ For a small neighborhood $U$ of $\alpha $ define $%
\phi :U\times F\rightarrow H(R^{2},R^{2})$ by $\phi (\beta ,f)=(h_{\alpha
\beta })h_{\alpha }f.$
\end{proof}

\begin{lemma}
\label{keylem}Suppose $[\alpha _{0}]=[id]\in \pi
_{1}(I(Z_{n},R^{2}),id_{|Z_{n}})$ and $\alpha _{t}$ is a path homotopy from $%
\alpha _{0}$ to $id$ in $I(Z_{n},R^{2}).$ Then there exists a level
preserving isotopy $h_{t}:R^{2}\times \lbrack 0,1]\rightarrow R^{2}\times
\lbrack 0,1]$ such that $h_{t|Z_{n}}=\alpha _{t},$ $h_{1}=id_{|R^{2}\times
\lbrack 0,1]},$ and $h_{t}$ is the identity on $R^{2}\times \{0\}$ and $%
R^{2}\times \{1\}$.
\end{lemma}

\begin{proof}
Define a map $q:[0,1]\times \lbrack 0,1]\rightarrow I(F,int(D))$ via $%
q(s,t)=\alpha _{t}(s).$ Note $id_{|F}=q(0,t)=q(s,1)=q(1,t).$ The space
obtained by collapsing 3 consecutive sides of $[0,1]\times \lbrack 0,1]$ to
a point is a topologically a disk. Thus by Proposition \ref{fib} $q$ lifts
to a map $Q:[0,1]\times \lbrack 0,1]\rightarrow H(R^{2},R^{2})$ such that $%
\pi (Q)=q$ and $id_{|D}=Q(0,t)=Q(s,1)=Q(1,t).$ Define $h_{t}:R^{2}\times
\lbrack 0,1]\rightarrow $ $R^{2}\times \lbrack 0,1]$ via $%
h_{t}(z,s)=Q(s,t)(z).$
\end{proof}

\begin{lemma}
\label{sigmamet}If $Z=\{1,2,3,...\}$ then $H(Z,Z)$ is a completely
metrizable topological group.
\end{lemma}

\begin{proof}
For $\{n,m\}\subset Z$ let $d^{\ast }(n,m)=1$ if $n\neq m$ and $d^{\ast
}(n,m)=0$ if $n=m.$ Note $C(Z,Z)$ is complete under the metric $%
d(f,g)=\sum_{n=1}^{\infty }\frac{d^{\ast }(f(n),g(n))}{2^{n}}.$ Let $%
E^{n}\subset C(Z,Z)$ denote the maps whose image contains $\{1,...n\}$.
Suppose $f\in E^{n}.$ Choose finite $M\subset Z$ such that $\{1,..n\}\subset
f(M)$ and $k=f(m_{k})$ with $1\leq k\leq n$ and $m_{k}\in M.$ Note for $g$
sufficiently close to $f$ we have $g(m_{k})=f(m_{k}).$ Hence $E^{n}$ is open
in $C(Z,Z)$. Thus $E(Z,Z)$ is a $G_{\delta }$ subspace of $C(Z,Z)$ since $%
E(Z,Z)=\cap E^{n}.$ In particular $E(Z,Z)$ is completely metrizable (p.270 
\cite{Munkres}). Let $I^{n}\subset C(Z,Z)$ denote the maps $f$ such that $%
f_{|\{1,..,n\}}$ is one to one. Suppose $f\notin I^{n}.$ Then for $g$
sufficiently close to $f$ we have $f_{|\{1,..,n\}}=g_{|\{1,..,n\}}.$ Hence $%
I^{n}$ is closed in $C(Z,Z).$ Thus $I(Z,Z)$ is closed in $C(Z,Z)$ since $%
I(Z,Z)=\cap I_{n}.$ Note $H(Z,Z)=I(Z,Z)\cap E(Z,Z).$ Hence $H(Z,Z)$ is
completely metrizable since it's a closed subspace of the completely
metrizable space $E(Z,Z).$

To prove $H(Z,Z)$ is a topological group note first that compact sets in $Z$
are exactly the finite subsets of $Z$. Hence convergence in $H(Z,Z)$ is
equivalent to pointwise convergence. If $f_{n}\rightarrow f$ and $%
g_{n}\rightarrow g$ then $\lim_{n\rightarrow \infty
}f_{n}(g_{n})(m)\rightarrow f(g)(m).$ Thus function composition is
continuous in $H(Z,Z).$ To check inversion is continuous it suffices to
check the case when $h_{n}\rightarrow id.$ Fixing $i$ choose $N$ such that $%
h_{n}(i)=i$ whenever $n\geq N.$ In particular $h_{n}^{-1}(i)=i.$ Thus $%
h_{n}^{-1}\rightarrow id,$ and hence $H(Z,Z)$ is a topological group.
\end{proof}

\begin{lemma}
$B_{\infty }$ is topologically incomplete with the topology described in the
introduction.
\end{lemma}

\begin{proof}
Note that $B_{\infty }$ is countable and that the complement of every one
point subset of $B_{\infty }$ is open and dense. Thus $B_{\infty }$ is the
countable union of closed nowhere dense sets. Hence by the Baire category
theorem $B_{\infty }$ does not admit a complete metric.
\end{proof}

\subsection{Generalizing the topological fundamental group}

It is unknown whether there exists a space whose topological fundamental
group completes $B_{\infty }$ with the topology described in the introduction%
$.$ To bypass this difficulty we develop a generalization of the topological
fundamental group. The basic idea is to replace the basepoint of $X$ with a
subset $A\subset X,$ and to consider paths of embeddings of $A$ into $X$
that initially fix $A$ pointwise, that are allowed to move $A$ off of itself
in midjourney, and that place $A$ back onto a permutation of itself at the
last moment. Under suitably nice conditions the path homotopy classes
determine a topological group.

Suppose throughout this section that $A\subset Y$, $\{id_{A}\}\subset
X\subset I(A,Y),$ and $H$ is a subgroup of $H(A,A).$

Let $K=\{f\in C([0,1],X)|f(0)=id_{A},$ $f(1)\in H\}.$

Define $m:K\times K\rightarrow K$ such that 
\[
m(f,g)(t)=\left\{ 
\begin{array}{ccc}
f(2t) &  & t\leq 1/2 \\ 
g(2t-1)f(1) &  & 1/2\leq 1
\end{array}
\right. 
\]

Define $v:K\rightarrow K$ such that $(vf)(t)=f(1-t)(f(1))^{-1}$ whenever $%
f\in K.$

Endowed with the quotient topology, let $G(X,H)$ denote the space whose
elements are path components of $K$.

\begin{theorem}
\label{topgrp}If $H$ is a topological group then $m$ and $v$ induce a
topological group structure on $G(X,H).$
\end{theorem}

\begin{proof}
Evidently $m$ is continuous and $v$ is continuous since inversion in $H$ is
continuous. Given $f\in K,$ let $[f]$ denote the path component of $f$ in $%
K. $ Let $q:K\rightarrow G(X,H)$ be the quotient map $q(f)=[f].$ Define $%
M:G(X,H)\times G(X,H)\rightarrow G(X,H)$ such that $M([f],[g])=[m(f,g)].$
Define $V:G(X,H)\rightarrow G(X,H)$ via $V([f])=[v(f)].$ To check that $M$
and $V$ are well defined and continuous it suffices by Proposition \ref{qmap}
to note that $q(m)$ and $q(j)$ are constant respectively on sets $(q\times
q)^{-1}([f],[g])$ and $q^{-1}([f]).$ Define $ID:[0,1]\rightarrow K$ via $%
ID(t)=id_{A}.$ The proof that $G(X,H)$ is a group with identity $[ID],$
multiplication $M,$ and inversion $V$ is similar to the proof that the
familiar fundamental group is indeed a group. Note the paths $%
m(ID,f),m(f,ID) $ and $f$ all differ by an order preserving surjective map
of $[0,1]$ and hence $M([ID],[f])=M([ID],[f])=[f].$ Similarly the maps $%
m(m(f,g),h)$ and $m(f,m(g,h))$ differ by an order preserving surjective map
of $[0,1]$ and thus $M$ is associative. Note $m(f,v(f))=m(v(f),f).$ The path 
$m(f,v(f))$ can be deformed to $ID$ by concatenating $f_{[s,0]}$ with $%
f_{[s,0]}.$ Thus $V$ serves as inversion and hence $G(X,H)$ is a topological
group.
\end{proof}

The based topological fundamental group $\pi _{1}(Y,p)$ \cite{Biss} is the
quotient space of loops based at $p$ from $[0,1]$ into $Y$ under the
equivalence relation path homotopy. This can be seen as a special case of
the generalized topological fundamental group $G(X,H).$

\begin{corollary}
If $A=\{p\}\subset Y,$ $H=\{id_{A}\}$ and $X=I(\{p\},Y)$ then $G(X,H)\cong
\pi _{1}(Y,p),$ the topological fundamental group of $X$ based at $p.$
\end{corollary}

\subsection{Twisted Products and Inverse Limit Spaces}

Lemma \ref{iso1} helps establish that the inverse limit of the groups $G_{n}$
is canonically isomorphic (and homeomorphic) to an iterated semidirect
product of groups $K_{n}.$

\begin{definition}
\label{twist}Suppose each of $G_{1},G_{2},..$ $K_{1},K_{2},..$ is a group
such that for all $n\geq 1$ $\cup _{i=1}^{n}K_{i}\subset G_{n}$ and $%
K_{n}\vartriangleleft G_{n}.$Then $\Pi _{n=1}^{\infty }K_{n}$ forms a group
under the binary operation determined coordinatewise by 
\[
\lbrack (x_{1},x_{2},..)\ast
(y_{1},y_{2},...)]_{n}=(y_{1}...y_{n-1})^{-1}(x_{n})(y_{1}..y_{n-1})y_{n} 
\]
(with multiplication performed in $G_{n})$. Call this group $\Pi
_{G_{n}}K_{n},$the \textbf{twisted product }of $K_{n}$ over $G_{n}$.
\end{definition}

\begin{lemma}
The twisted product is a well defined group.
\end{lemma}

\begin{proof}
Note $\ast $ is well defined since $K_{n}\vartriangleleft G_{n}$ and $%
K_{i}\subset G_{n}$ for $1\leq i<n.$ Checking associativity, 
\[
\lbrack x\ast (y\ast
z)]_{n}=[(y_{1}z_{1})(z_{1}^{-1}y_{2}z_{1}z_{2})((z_{1}z_{2})^{-1}y_{3}(z_{1}z_{2}z_{3}))... 
\]
\[
((z_{1}...z_{n-2})^{-1}y_{n-1}z_{1}z_{2}...z_{n-1})]^{-1}x_{n}[(y_{1}z_{1})(z_{1}^{-1}y_{2}z_{1}z_{2})((z_{1}z_{2})^{-1}y_{3}(z_{1}z_{2}z_{3}))... 
\]
\[
((z_{1}...z_{n-2})^{-1}y_{n-1}z_{1}z_{2}...z_{n-1})](z_{1}...z_{n-1})^{-1}y_{n}(z_{1}...z_{n})= 
\]
\[
(z_{1}..z_{n-1})^{-1}(y_{1}...y_{n-1})^{-1}x_{n}(y_{1}..y_{n})(z_{1}..z_{n})=[(x\ast y)\ast z]_{n}. 
\]

Let $ID=(id,id,..).$ Then $[IDx]_{n}=x_{n}=[xID]_{n}.$ Finally given $%
y=(y_{1},y_{2},..)$ let $x=(x_{1},x_{2},..)$ with $%
x_{n}=(y_{1}...y_{n-1})(y_{1}...y_{n-1}y_{n})^{-1}.$ Then 
\[
\lbrack x\ast
y]_{n}=(y_{1}..y_{n-1})^{-1}(y_{1}...y_{n-1})(y_{n})^{-1}(y_{1}...y_{n-1})^{-1}(y_{1}...y_{n})=id. 
\]
\end{proof}

\begin{lemma}
\label{twistcor}Suppose $G$ is a group with subgroups $\{id\}=G_{0}\subset
G_{1}\subset G_{2}\subset G_{3}...\subset G$ and there exist homomorphisms $%
\phi _{n}:G\rightarrow G_{n}$ for all $n\geq 0$ satisfying $\phi
_{n|G_{n}}=id_{|G_{n}}.$ Define for all $n\geq 1$ $\psi
_{n}:G_{n}\rightarrow G_{n-1}$ as $\psi _{n}=\phi _{n-1|G_{n}}.$ Let $%
K_{n}=\ker \psi _{n}.$ Then $\phi :G\rightarrow \Pi _{G_{n}}K_{i}$ defined
as $(\phi (g))_{n}=(\phi _{n-1}(g))^{-1}\phi _{n}(g)$ is a homomorphism.
\end{lemma}

\begin{proof}
Note $K_{n}=\ker \psi _{n}\vartriangleleft G_{n}$ and $K_{i}\subset
G_{i}\subset G_{n}$ for $1\leq i<n$ and the hypothesis of Definition \ref
{twist} is satisfied. To check $im(\phi )\subset \Pi K_{i}$ we verify $(\phi
_{n-1}(g))^{-1}\phi _{n}(g)\in K_{n}.$ Let $p=\phi _{n}(g)$. Note $\phi
_{n-1}(g)=\phi _{n-1}(\phi _{n}(g))=\phi _{n-1}(p)$ and $\phi _{n-1}(\phi
_{n-1}(x))=\phi _{n-1}(x)$ for all $x\in G.$ Thus 
\[
\psi _{n}((\phi _{n-1}(p))^{-1}p)=\phi _{n-1}(\phi _{n-1}(p^{-1}))\phi
_{n-1}(p)=(\phi _{n-1}(p^{-1}))\phi _{n-1}(p)=id. 
\]
To check $\phi $ is a homomorphism use the following notation. For $g\in G$
let $g_{n}=(\phi (g))_{n}\in K_{n}$ and let $g^{(n)}=\phi _{n}(g)\in G_{n}.$
Note 
\[
g^{(n)}=(g^{(0)})^{-1}(g^{(1)})(g^{(1)})^{-1}..(g^{(n-1)})^{-1}g^{(n)}=g_{1}g_{2}...g_{n}. 
\]
Thus 
\[
(\phi (gh))_{n}=(\phi _{n-1}(gh))^{-1}\phi _{n}(gh)=\phi _{n-1}(h^{-1})(\phi
_{n-1}(g^{-1})\phi _{n}(g)\phi _{n}(h)= 
\]
\[
(h^{(n-1)})^{-1}g_{n}h^{(n)}=(h_{1}..h_{n-1})g_{n}(h_{1}...h_{n-1})h_{n}=(%
\phi (g)\ast \phi (h))_{n}. 
\]
\end{proof}

\begin{lemma}
\label{iso1}Suppose $\{id\}=G_{0}\subset G_{1}\subset G_{2}..$is a sequence
of groups and $\psi _{n}:G_{n}\rightarrow G_{n-1}$ is an epimorphism such
that $\psi _{n|G_{n-1}}=id_{|G_{n-1}}.$ Let $\lim_{\leftarrow
}G_{n}=\{(g_{1},g_{2},..)|g_{n}\in G_{n}$ and $\psi _{n}(g_{n})=g_{n-1}\}$
with group operation $(g_{1},g_{2},..)\ast
(h_{1},h_{2},..)=(g_{1}h_{1},g_{2}h_{2},..)$. Let $K_{n}=\ker \psi _{n}.$
Then $\psi :\Pi _{G_{n}}K_{n}\rightarrow \lim_{\leftarrow }G_{n}$ defined as 
$\psi (k_{1},k_{2},..)=(k_{1},k_{1}k_{2},..)$ is an isomorphism. If for all $%
n$ $G_{n}$ has the discrete topology and $\lim_{\leftarrow }G_{n}$ and $\Pi
_{G_{n}}K_{n}$ have the product topology then $\psi $ is a homeomorphism.
\end{lemma}

\begin{proof}
To check $\psi $ is a homomorphism 
\[
\psi \lbrack (x_{1},x_{2},...)(y_{1},y_{2},..)]=\phi
(x_{1}y_{1},(y_{1})^{-1}x_{2}y_{1}y_{2},...)= 
\]
\[
x_{1}y_{1},x_{1}x_{2}y_{1}y_{2},...=(x_{1},x_{1}x_{2},..)\ast
(y_{1},y_{1}y_{2},..)=\psi (x_{1},x_{2},..)\ast \psi (y_{1},y_{2},..). 
\]

Suppose $\psi (k_{1},k_{2},..)=(id,id,..).$ Then $k_{1}=id,$ $%
k_{1}k_{2}=id,...$ Hence by induction $k_{n}=id$ for all $n.$ Suppose $%
(g_{1},g_{2},..)\in \lim_{\leftarrow }G_{n}.$ Let $k_{1}=g_{1}.$ Suppose by
induction that $g_{n-1}=k_{1}...k_{n-1}$ with $k_{i}\in K_{i}.$ Let $%
k_{n}=(g_{n-1})^{-1}g_{n}.$

Note 
\[
\psi _{n}(k_{n})=\psi _{n}((g_{n-1})^{-1}\psi _{n}(g_{n})=g_{n-1}^{-1}\psi
(g_{n})=\psi (g_{n}^{-1})\psi (g_{n})=id. 
\]
Thus $\psi $ is an isomorphism. Note $\lim_{\leftarrow }G_{n}$ and $\Pi
_{G_{n}}K_{n}$ are metrizable since each of $G_{n}$ and $K_{n}$ is
metrizable. Thus it suffices to prove the homomorphisms $\psi $ and $\psi
^{-1}$ preserve convergence at $(id,id,..).$ Suppose 
\[
(k_{1}^{1},k_{2}^{1},k_{3}^{1},...),(k_{1}^{2},k_{2}^{2},k_{3}^{2},..),...%
\rightarrow (id,id,id,..). 
\]
Suppose $N>0.$ Choose $M$ so that $k_{j}^{i}=id$ for $i\geq M$ and $j\leq N.$
Suppose $i\geq M$ then $[\psi (k_{1}^{i},k_{2}^{i},..)]_{n}=id.$ Thus $\psi $
is continuous. Conversely suppose 
\[
(g_{1}^{1},g_{2}^{1},g_{3}^{1},..),(g_{1}^{2},g_{2}^{2},g_{3}^{2},..),...%
\rightarrow (id,id,id,...). 
\]
Suppose $N>0.$ Choose $M$ so that $i\geq M$ and $g_{N}^{i}=g_{N-1}^{i}=id.$
Then $k_{n}=id.$ Thus $\psi ^{-1}$ is continuous.
\end{proof}

\section{\label{gcomplete}The complete braid groups $G$ and $\overline{%
B_{\infty }}$ on infinitely many strands}

Let $X=I(Z^{+},R^{2})$

Let $X_{n}=\{f\in X|$ $f(i,0)=(i,0)$ \ for all $i>n$ and for all $i\leq n$ $%
f(i,0)\in (-\infty ,n+2/3]\times R\}.$

Recalling Theorem \ref{topgrp}, let $G=G(X,\{id_{Z^{+}}\}),$ the topological
fundamental group of $X$ based at $id_{Z^{+}}.$

Viewed as a collection of arcs in $R^{2}\times \lbrack 0,1],$ given $g\in
\lbrack g]\in G,$ the strands of $g$ connect $(n,0,0)$ to $(n,0,1)$ for $%
n\in \{1,2,3,...\}.$

The main goal is to prove that $G$ is isomorphic to the inverse $%
\lim_{\leftarrow }G_{n}$ where $G_{n}$ is the pure braid group on $n$
strands, realized as the image of $G$ under retractions that move the first $%
n$ strands `out of the way' and then replace the remaining strands with
trivial strands.

\subsection{The retraction $P_{n}:X\rightarrow X_{n}$}

\begin{lemma}
\label{retract}There exists a retraction $P_{n}:X\rightarrow X_{n}.$
\end{lemma}

\begin{proof}
Define $R_{n}:X\rightarrow \lbrack 0,\infty )$ as $R_{n}(g)=\min \{r|r\geq 0$
and for all $i\leq n$ $g(i,0)\in (-\infty ,n+2/3+r]\times R\}.$ Define $%
T:R\rightarrow C(R^{2},R^{2})$ as $T(r)(x,y)=(x-r,y).$ Define $%
P_{n}:X\rightarrow X_{n}$ via $P_{n}(g)(i,0)=\left\{ 
\begin{array}{ccc}
(TR_{n}g)(g(i,0)) &  & i\leq n \\ 
(i,0) &  & i>n
\end{array}
\right. .$
\end{proof}

\subsection{\label{defgn}The groups $G_{n}$ and homomorphisms $\protect\phi %
_{n}$ and $\protect\phi $}

Suppose $n\geq 1.$ Define $P_{n}^{\ast }:\pi _{1}(X,id)\rightarrow \pi
_{1}(X_{n},id)$ as $P_{n}^{\ast }[g]=[P_{n}(g)]$ where $P_{n}:X\rightarrow
X_{n}$ is the retraction from Lemma \ref{retract}. Let $i_{n}^{\ast }:\pi
_{1}(X_{n},id)\rightarrow \pi _{1}(X,id)$ satisfy $i_{n}^{\ast }[\alpha
]=[\alpha ].$ Recall $G=\pi _{1}(X,id_{Z^{+}}).$ Let $G_{n}=im(i_{n}^{\ast
})\subset G.$ Define $\phi _{n}:G\rightarrow G_{n}$ as $\phi
_{n}=i_{n}^{\ast }P_{n}^{\ast }.$ Let $G_{0}=\{[id]\}\subset G$ and define $%
\phi _{0}:G\rightarrow G_{0}$. Let $K_{n}=G_{n}\cap \ker \phi _{n-1}.$
Define $\phi :G\rightarrow \Pi _{G_{n}}K_{n}$ as in Lemma \ref{twistcor}.

\subsection{$\protect\phi :G\rightarrow \Pi _{G_{n}}K_{n}$ is onto}

\begin{lemma}
\label{nnoti}Suppose $n\geq 1$ and $[k_{n}]\in K_{n}.$ Then there exists $%
\beta \in \lbrack k_{n}]$ such that $\beta (s)(i,0)=(i,0)$ for all $s\in
\lbrack 0,1]$ and for all $i\leq n-1.$
\end{lemma}

\begin{proof}
Since $[k_{n}]\in im(i_{n}^{\ast }),$ we may choose $\alpha \in \lbrack
k_{n}]$ such that the $ith$ strand of $\alpha $ is straight for all $i>n.$
The braid $\alpha $ can deformed in $X_{n}$ to a braid $\Gamma $ so that the
first $n-1$ strands of $\Gamma $ lie over the closed half plane $(\infty ,n-%
\frac{1}{3}]\times R$ as follows. Recalling the map $R_{n}$ from the proof
of Lemma \ref{retract}, fix $t$ and translate $\alpha (t)(Z_{n})$ leftward
at constant speed by a total amount $R_{n-1}(\alpha (t)).$ The braid $\gamma
=P_{n-1}(\Gamma )$ is equivalent to the trivial braid in $X_{n-1}$ and the
first $n-1$ strands of $\gamma $ agree with the first $n-1$ strands of $%
\Gamma .$ Let $\gamma _{t}$ be a path of braids in $X_{n-1}$ taking $\gamma $
to $id.$ By Lemma \ref{keylem} extend $\gamma _{t}$ to a level preserving
isotopy $h_{t}:R^{2}\times \lbrack 0,1]\rightarrow R^{2}\times \lbrack 0,1]$
such that and $h_{1}=id_{|R^{2}\times \lbrack 0,1]}$ and $h_{t}$ is the
identity on the top and bottom. Now consider the path of braids $\Gamma
_{t}=h_{t}h_{0}^{-1}(\Gamma ).$ Note $\Gamma _{0}=\Gamma $ and the first $%
n-1 $ strands of $\Gamma _{1}$ are straight. Let $\beta =\Gamma _{1}.$
\end{proof}

\begin{corollary}
\label{epi}$\phi $ is an epimorphism.
\end{corollary}

\begin{proof}
Suppose $\{[k_{1}],[k_{2}],..\}\in \Pi _{G_{n}}K_{i}.$

Choose by Lemma \ref{nnoti} $\beta _{n}\in \lbrack k_{n}]$ such that $\beta
_{n}(s)(i,0)=(i,0)$ for $i\leq n-1.$ Define $g:[1,\infty ]\rightarrow X$ as

$g(n+t)=\left\{ 
\begin{array}{ccc}
\beta _{n}(t) &  & n\geq 1\text{ and }t\in \lbrack 0,1] \\ 
\iota d &  & t=\infty
\end{array}
\right. .$

Note $g$ is well defined and continuous when $t<\infty $ since $\beta
_{n}(0)=\beta _{n}(1)=id$ for $n\geq 1.$ Suppose $t=\infty $, $%
T_{n}\rightarrow \infty $ and $F\subset Z^{+}$ is compact. Let $M=\max
\{i|(i,0)\in F\}.$ For $T\geq M+1$ and $(i,0)\in F,$ $g(T)(i,0)=(i,0).$
Hence $g$ is continuous. Let $h:[0,1]\rightarrow \lbrack 1,\infty ]$ be any
order preserving homeomorphism. Let $f=gh.$ Note $f(0)=f(1)=id.$ Note $\phi
_{n}([f])=[\beta _{1}\beta _{2}...\beta _{n}].$ Thus, recalling Lemma \ref
{twistcor}, $(\phi ([f])_{n}=\phi _{n-1}([f])^{-1}\phi _{n}([f])=[\beta
_{1}..\beta _{n-1}]^{-1}[\beta _{1}\beta _{2}...\beta _{n}]=[\beta _{n}].$
Hence $\phi $ is an epimorphism.
\end{proof}

\subsection{$\protect\phi :G\rightarrow \Pi _{G}K_{n}$ is one to one}

Lemma \ref{kill1} says the following: Given an infinite pure braid $\alpha ,$
if the infinite braid $P_{n+1}(\alpha )$ is equivalent to the trivial braid
on infinitely many strands, then $\alpha _{|Z_{n+1}}$ is equivalent to the
trivial braid on $n+1$ strands.

\begin{lemma}
\label{kill1}Suppose $[\alpha ]\in \ker \phi _{n+1}.$ Define $\beta
:[0,1]\rightarrow I(Z_{n+1},R^{2})$\ as $\beta (s)=\alpha (s)_{|Z_{n+1}}.$
Then $[\beta ]=[id]\in \pi _{1}(I(Z_{n+1},R^{2}),id).$
\end{lemma}

\begin{proof}
For $(s,t)\in \lbrack 0,1]\times \lbrack 0,1]$ define $\beta
_{s}^{t}:Z_{n+1}\hookrightarrow R^{2}$ via

$\beta ^{t}(s)(i,0)=(T(tR_{n+1}(\alpha (s))))\alpha (i,0).$ Note $\beta
^{0}=\beta ,P_{n+1}(\alpha )_{|Z_{n+1}}=\beta ^{1}$ and $\beta ^{1}\in
\lbrack \beta ]\in \pi _{1}(I(Z_{n+1},R^{2}),id).$ Since $[P_{n+1}(\alpha
)]=[id]\in \pi _{1}(X_{n},id)$ there exists $\gamma :[0,1]\rightarrow X_{n}$
such that $\gamma (0)=P_{n+1}(\alpha )$ and $\gamma (1)=id.$ Now observe $%
\gamma (t)_{|Z_{n+1}}$ connects $\beta ^{1}$ to $id$ within $[\beta ]\in \pi
_{1}(I(Z_{n+1},R^{2}),id).$
\end{proof}

Lemma \ref{kill2} says the following. Suppose the pure braid on $n+1$
strands $\beta ^{0}$ can be deformed to the trivial braid on $n+1$ strands,
and $n$ strands of $\beta ^{0}$ are straight. Then the deformation can be
chosen to leave these $n$ strands invariant.

\begin{lemma}
\label{kill2}Suppose $[\beta ^{0}]=[id]\in \pi _{1}(I(Z_{n+1},R^{2}),id)$
and for all $s\in \lbrack 0,1]$ $\beta ^{0}(s)_{Z_{n}}=id_{Z_{n}}.$ Then
there exists $\gamma :[0,1]\rightarrow \lbrack \beta ^{0}]$ such that $%
\gamma (0)=\beta ^{0},\gamma (1)=id$ and $\gamma (t)_{Z_{n}}=id_{Z_{n}}.$
\end{lemma}

\begin{proof}
Let $f:[0,1]\times \lbrack 0,1]\rightarrow I(Z_{n+1},R^{2})$ satisfy $%
f(0,s)=\beta ^{0}(s)$ and $f(1,s)=f(0,t)=f(1,t)=id.$ Let $%
q(s,t)=f(s,t)_{Z_{n}}.$ Note for $(s,t)\in \partial ([0,1]\times \lbrack
0,1])$ $q(s,t)_{Z_{n}}=id_{Z_{n}}.$ Since $I(Z_{n},R^{2})$ is aspherical
(Theorem 2\cite{Fadell}), viewed as a map of $S^{2}$ into $I(Z_{n},R^{2}),$ $%
q$ is inessential. By Proposition \ref{fib} the restriction map $\pi
:H(R^{2},R^{2})\rightarrow I(Z_{n},R^{2})$ defined as $\pi (h)=h_{Z_{n}}$
determines a fibre bundle. Hence there exists a `lift' $Q:[0,1]\times
\lbrack 0,1]\rightarrow H(R^{2},R^{2})$ such that $Q(s,t)=ID_{R^{2}}$ for $%
(s,t)\in \partial ([0,1]\times \lbrack 0,1])$ and $\pi (Q)=q.$ Let $\gamma
(t)(s)=(Q(s,t))^{-1}f(s,t).$
\end{proof}

Lemma \ref{kill3} establishes the inductive step of our proof that $\phi $
is one to one.

\begin{lemma}
\label{kill3}Suppose $[\alpha _{0}]\in \ker \phi _{n+1}$ and for all $s\in
\lbrack 0,1]$ $\alpha _{0}(s)_{|Z_{n}}=id_{Z_{n}}.$ Then there exists a map $%
\Gamma :[0,1]\rightarrow \lbrack \alpha _{0}]$ such that $\Gamma (0)=\alpha
_{0},$ $\Gamma (t)(s)_{|Z_{n}}=id_{|Z_{n}}$ and $\Gamma
(1)(s)_{|Z_{n+1}}=id_{|Z_{n+1}}.$
\end{lemma}

\begin{proof}
Define $\beta _{0}:[0,1]\rightarrow I(Z_{n+1},R^{2})$ via $\beta
_{0}(s)=\alpha _{0}(s)$ applying Lemmas \ref{kill1} and \ref{kill2} let $%
\beta _{t}$ be a path homotopy from $\beta _{0}$ to $id$ in $%
I(Z_{n+1},R^{2}) $ such that the first $n$ strands of $\beta _{t}$ are
straight. By Lemma \ref{keylem} extend $\beta _{t}$ to a level preserving
isotopy $h_{t}:R^{2}\times \lbrack 0,1]\rightarrow R^{2}\times \lbrack 0,1]$
such that $h_{1}=id_{|R^{2}x[0,1]}$ and $h_{t}$ is the identity on the top
and bottom level. Let $\Gamma _{t}=h_{t}h_{0}^{-1}(\alpha _{0}).$
\end{proof}

\begin{lemma}
$\phi :G\rightarrow \Pi _{G}K_{n}$ is one to one
\end{lemma}

\begin{proof}
By definition $\ker (\phi )=\cap _{i=n}^{\infty }\ker (\phi _{n}).$ To prove 
$\phi $ has trivial kernel we take $[g_{0}]\in \ker \phi $ and build a path $%
g_{t}$ (parametrized over $[0,\infty ]$) from $g_{0}$ to $id$ within $[g_{0}]
$ by the following procedure. By induction assume we have defined the
homotopy $g_{t}$ for $t\in \lbrack 0,n]$ and that the first $n$ strands of $%
g_{n}$ are straight. Then straighten the $n+1$ strand of $g_{n}$ for $t\in
\lbrack n,n+1]$ via methods established in Lemmas \ref{kill1}, \ref{kill2},
and \ref{kill3}. Define $g_{\infty }(s)=id_{|Z^{+}}.$ We must prove the
function $q:[0,\infty ]\rightarrow C([0,1],C(Z^{+},R^{2}))$ defined via $%
q(t)=g_{t}$ is continuous when $t=\infty .$ By Proposition \ref{munkcor} it
suffices to show that the following map $Q:[0,\infty ]\times \lbrack
0,1]\times Z^{+}\rightarrow R^{2}$ defined via $Q(t,s,(n,0))=g_{t}(s)(n,0)$
is continuous. Continuity is immediate if $t<\infty .$ Suppose $%
t_{k}\rightarrow \infty $ and $s_{k}\rightarrow s$ and $(n_{k},0)\rightarrow
(n,0).$ Suppose $\varepsilon >0.$ Choose $K$ such that if $k\geq K$ then $%
n_{k}=n$ and $t_{k}>n+1.$ Suppose $k\geq K.$ Then $%
Q(t_{k},s_{k},n_{k})=(n,0)=Q(\infty ,s,n).$ Thus $q$ is continuous.
\end{proof}

\subsection{$G$ is the inverse limit of the pure braid groups}

\begin{theorem}
\label{gcomp}$\phi :G\rightarrow \Pi _{G_{n}}K_{n}$ is a homeomorphism and $%
G $ is topologically complete.
\end{theorem}

\begin{proof}
To prove $\phi _{n}:G\rightarrow G_{n}$ is continuous apply Proposition \ref
{qmap} to the retraction $P_{n}:X\rightarrow X_{n}.$ Thus the isomorphism $%
H^{\ast }:G\rightarrow \lim_{\leftarrow }G_{n}$ defined via $H^{\ast
}([g])=(\phi _{1}([g]),\phi _{2}([g]),...)$ is continuous. Note $\phi =\psi
H^{\ast }$ where $\psi :\lim_{\leftarrow }G_{n}\rightarrow \Pi _{G_{n}}K_{n}$
is the homeomorphism from Lemma \ref{iso1}. Thus $\phi $ is continuous. To
check $\phi ^{-1}$ is continuous it suffices, since $\Pi _{G_{n}}K_{n}$ is
metrizable and $\phi ^{-1}$ is a homomorphism, to show $\phi ^{-1}$
preserves convergence at $([id],[id],...).$ Suppose $%
([k_{1}^{1}],[k_{2}^{1}],...),([k_{1}^{2}],[k_{2}^{2}],..),...\rightarrow
([id],[id],..)\in \Pi _{G_{n}}K_{n}.$ Choose $g^{i}\in \phi
^{-1}([k_{1}^{i}],[k_{2}^{i}],...)$ such that ( considering $im(g^{i})$ as a
subspace of $R^{2}\times \lbrack 0,1])$ a maximal number of initial strands
of $im(g)$ are straight line segments. It follows that $g^{i}\rightarrow
id^{\ast }$ in $C([0,1],X)$ where $id^{\ast }(t)=id_{Z^{+}}.$

Hence $[g^{i}]\rightarrow \lbrack id^{\ast }]\in G$ and thus $\phi ^{-1}$ is
continuous. Moreover $G$ is completely metrizable since $G$ is homeomorphic
to the countable product of complete (discrete) spaces $K_{n}.$
\end{proof}

\subsection{The complete braid group $\overline{B_{\infty }}$}

Recalling Theorem \ref{topgrp} let $\overline{B_{\infty }}%
=G(X,H(Z^{+},Z^{+})).$

In the definition of $G(X,H(Z^{+},Z^{+})),$ when defining $K,$ for
convenience replace $[0,1]$ with $[0,\infty ]$ the one point
compactification of $[0,\infty )$.

Elements of $\overline{B_{\infty }}$ are path components of $K,$ and the
graph of each element of $K$ can be seen as a collection of arcs in $%
R^{2}\times \lbrack 0,\infty ]$ connecting $Z^{+}\times \{0\}$ to $%
Z^{+}\times \{\infty \}.$ Roughly speaking, two such graphs are equivalent
if one can be deformed into the other while leaving the endpoints fixed.

We show $\overline{B_{\infty }}$ is topologically complete by exhibiting a
homeomorphism onto the product $G\times H(Z^{+},Z^{+}).$

\begin{lemma}
\label{sect}Define $\Pi :\overline{B_{\infty }}\rightarrow H(Z^{+},Z^{+})$
such that $\Pi ([b])(i)=b(\infty )(i).$ Suppose there exists a map $\sigma
:H(Z^{+},Z^{+})\rightarrow \overline{B_{\infty }}$ such that $\Pi \sigma
=id. $ Then $\overline{B_{\infty }}$ is homeomorphic to the product $%
G(X,id_{Z^{+}})\times H(Z^{+},Z^{+}).$
\end{lemma}

\begin{proof}
Define $T:K\rightarrow H(Z^{+},Z^{+})$ via $T(k)=k(\infty ).$ Note all path
components of $H(Z^{+},Z^{+})$ are trivial and thus $T$ induces the map $\Pi 
$ by Proposition \ref{qmap}. Moreover $\Pi $ is a homomorphism and $\ker \Pi 
$ is precisely $G(X,id_{Z^{+}}),$ the fundamental group of $X$ based at $%
id_{Z^{+}}.$ Now consider the map $F:\ker (\Pi )\times
H(Z^{+},Z^{+})\rightarrow \overline{B_{\infty }}$ defined via $F(g,\tau
)=g\ast \sigma (\tau ).$ Consider the map $P:\overline{B_{\infty }}%
\rightarrow \ker (\Pi )\times H(Z^{+},Z^{+})$ defined via $P(b)=(b\ast
(\sigma \Pi b)^{-1},\Pi b).$

To check that $b\ast (\sigma \Pi b)^{-1}\in \ker (\Pi )$ note $\Pi (b\ast
(\sigma \Pi b)^{-1})=\Pi (b)\ast \Pi ((\sigma \Pi b)^{-1})=\Pi (b)\ast
\lbrack \Pi (\sigma \Pi b)]^{-1}=\Pi (b)\ast \lbrack \Pi b]^{-1}=id.$
Finally note $FP=id.$ Thus each of $F$ and $P$ is a homeomorphism.
\end{proof}

Let $Z=\{1,2,3,...\}$ with the discrete topology$.$ We will show that $%
H(Z,Z) $ is a closed subspace of $E(Z,Z)\subset C(Z,Z)$, that $E(Z,Z)$ is a $%
G_{\delta }$ subspace of $C(Z,Z),$ and that $C(Z,Z)$ admits a complete
metric. It follows that $H(Z,Z)$ admits a complete metric.

We now construct a `normal form' for each permutation $\tau \in H(Z,Z)$ as
an `infinite word' $\sigma _{1m_{1}}\sigma _{2m_{2}}...$ where $\sigma
_{nm}\in \sigma _{\infty }$ has compact support.

For $1\leq n\leq m\in Z^{+}$ define $\sigma _{nm}:Z^{+}\rightarrow Z^{+}$via 
\[
\sigma _{nm}(i)=\left\{ 
\begin{array}{ccc}
i+1 &  & n\leq i\leq m-1 \\ 
n &  & i=m \\ 
i &  & \text{otherwise}
\end{array}
\right. 
\]

Define $b_{nm}:[0,1]\rightarrow X$ such that $b_{nm}$ is continuous and

\[
\left\{ 
\begin{array}{ccc}
b_{nm}(1)(i,0)=\sigma _{nm}(i) &  & \text{ for all }i \\ 
b_{nm}(0)(i,0)=(i,0) &  & \text{ for all }i \\ 
n-1<\left| b_{nm}(t)(i,0)\right| <m+1 &  & n\leq i\leq m\text{ and }n<m \\ 
b_{nm}(t)(i,0)=(i,0) &  & \text{otherwise}
\end{array}
\right. 
\]

\begin{lemma}
\label{normalperm}For each $\tau \in H(Z,Z)$ there exists a sequence 
\[
\sigma _{1m_{1}},\sigma _{2m_{2}},..\in H(Z,Z) 
\]
such that for all $k\geq 1$ 
\[
(\tau ^{-1})_{|\{1,...k\}}=(\sigma _{km_{k}}...\sigma _{2m_{2}}\sigma
_{1m_{1}})_{|\{1,..k\}}^{-1}. 
\]
\end{lemma}

\begin{proof}
Suppose $\tau \in H(Z,Z).$ Let $m_{1}=\tau ^{-1}(1).$ Proceeding by
induction suppose

$(\tau ^{-1})_{|\{1,...k\}}=(\sigma _{km_{k}}...\sigma _{2m_{2}}\sigma
_{1m_{1}})_{|\{1,..k\}}^{-1}$.

Define $m_{k+1}=\sigma _{km_{k}}...\sigma _{2m_{2}}\sigma _{1m_{1}}\tau
^{-1}(k+1).$ Since $\sigma _{km_{k}}...\sigma _{2m_{2}}\sigma _{1m_{1}}$
maps $\tau ^{-1}\{1,.,k\}$ onto $\{1,..k\}$ it follows that $m_{k+1}\geq
k+1. $ By the induction hypothesis if $i\leq k$ then $(\sigma
_{k+1m_{k+1}}\sigma _{km_{k}}...\sigma _{2m_{2}}\sigma
_{1m_{1}})^{-1}(i)=\sigma _{1m_{1}}^{-1}..\sigma _{km_{k}}^{-1}(i)=\tau
^{-1}(i)$ since $\sigma _{k+1m_{k+1}}$ fixes $\{1,..k\}$ pointwise. If $%
i=k+1 $ then 
\[
(\sigma _{k+1m_{k+1}}\sigma _{km_{k}}...\sigma _{2m_{2}}\sigma
_{1m_{1}})^{-1}(k+1)= 
\]
\[
\sigma _{1m_{1}}^{-1}..\sigma _{km_{k}}^{-1}\sigma
_{k+1m_{k+1}}^{-1}(k+1)=\sigma _{1m_{1}}^{-1}..\sigma
_{km_{k}}^{-1}(m_{k+1})= 
\]
\[
\sigma _{1m_{1}}^{-1}..\sigma _{km_{k}}^{-1}(\sigma _{km_{k}}...\sigma
_{2m_{2}}\sigma _{1m_{1}}\tau ^{-1}(k+1))=\tau ^{-1}(k+1). 
\]
\end{proof}

\begin{lemma}
\label{section}There exists a continuous function $\sigma :H(Z,Z)\rightarrow 
\overline{B_{\infty }}$ such that $\sigma (\tau )(1)(i,0)=(\tau (i),0)$ for
all $i\in H(Z,Z).$
\end{lemma}

\begin{proof}
Define $\stackrel{\symbol{94}}{\sigma }:H(Z,Z)\rightarrow K\subset
C([0,\infty ],X)$ as follows.

Given $\tau \in H(Z,Z)$ choose $\sigma _{1m_{1}},\sigma _{2m_{2}},..$ via
Lemma \ref{normalperm} so that for all $k\geq 1(\tau
^{-1})_{|\{1,...k\}}=(\sigma _{km_{k}}...\sigma _{2m_{2}}\sigma
_{1m_{1}})_{|\{1,..k\}}^{-1}.$ Define $\alpha _{\tau }:[0,\infty
]\rightarrow X$ such that 
\[
\alpha _{\tau }(t)(i,0)=\left\{ 
\begin{array}{ccc}
(i,0) &  & t=0 \\ 
b_{km_{k}}(t)\alpha _{\tau }(k-1)(i,0) &  & t\in (k-1,k] \\ 
(\tau (i),0) &  & t=\infty
\end{array}
\right. 
\]
Note $\alpha _{\tau }$ is continuous at $k$ since $b_{km_{k}}(0)=id.$ To
check continuity at $t=\infty .$ Suppose $(i,0)\in Z^{+}.$ Let $(n,0)=(\tau
(i),0).$ Suppose $t\in (k,k+1]$ and $n+1\leq k.$ Then $\alpha _{\tau
}(k+t)(i,0)=b_{km_{k}}(t)b_{k-1m_{k-1}}...b_{1m_{1}}(1)(i,0)=b_{km_{k}}(t)(%
\sigma _{k-1m_{k-1}}...\sigma _{1m_{1}}(i),0)=b_{km_{k}}(t)(n,0)=(n,0).$
Thus $\alpha _{\tau }$ is continuous on $[0,\infty ].$

To check $\stackrel{\symbol{94}}{\sigma }$ is continuous suppose $\tau
_{m}\rightarrow \tau \in H(Z,Z).$

Choose $\sigma _{1m_{1}},\sigma _{2m_{2}},..$ via Lemma \ref{normalperm} so
that for all $k\geq 1(\tau ^{-1})_{|\{1,...k\}}=(\sigma _{km_{k}}...\sigma
_{2m_{2}}\sigma _{1m_{1}})_{|\{1,..k\}}^{-1}.$ Suppose $(i,0)\in Z^{+}.$ Let 
$n=\tau (i).$ Choose $M>n$ so that if $m\geq M$ then $\sigma
_{km_{k}}^{m}=\sigma _{km_{k}}$ where $\sigma _{1m_{1}}^{m},\sigma
_{2m_{2}}^{m}...$ is chosen for $\tau _{m}$ as in Lemma \ref{normalperm}.
Suppose $m\geq M$ and $t\in \lbrack 0,\infty ].$ Note $\alpha _{\tau
_{m}}(t)(i,0)=\alpha _{\tau }(t)(i,0).$ Thus $\sigma ^{\symbol{94}}$ is
continuous.

Appealing to Proposition \ref{qmap} let $\sigma :\Sigma _{\infty
}\rightarrow \overline{B_{\infty }}$ be the induced map defined via $\sigma
(\tau )=[\stackrel{\symbol{94}}{\sigma }(\tau )].$
\end{proof}

\begin{theorem}
\label{topcomp}$\overline{B_{\infty }}$ is topologically complete.
\end{theorem}

\begin{proof}
Applying Lemmas \ref{sect}, \ref{sigmamet}, \ref{section} and Theorem \ref
{gcomp} we obtain a homeomorphism between $\overline{B_{\infty }}$ and the
product of topologically complete spaces $G(X,id_{Z^{+}})\times
H(Z^{+},Z^{+}).$
\end{proof}

\subsection{A geometric description of convergence in $\overline{B_{\infty }}%
.$}

We indicate why convergence $b_{n}\rightarrow id\in \overline{B_{\infty }}$
is characterized by the statement `For each $M$ there exists $N$ such that
if $n\geq N$ then the first $M$ strands of $b_{n}$ are trivial'. Recalling
the proof of Lemma \ref{sect} note $b_{n}\rightarrow id$ iff $\Pi
(b_{n})\rightarrow id\in H(Z^{+},Z^{+})$ and $b_{n}\ast (\sigma \Pi
b_{n})^{-1}\rightarrow id\in G.$ In $H(Z^{+},Z^{+})$ convergence is
pointwise. Thus $\Pi (b_{n})\rightarrow id$ iff for each $M$ there exists $N$
such that $n\geq N$ $\Pi (b_{n})(i)=i$ for $i\leq M$ iff (recalling the
definition of $\sigma $ from Lemma \ref{section}) for each $M$ there exists $%
N$ such that the first $M$ strands of $\sigma (\Pi b_{n})$ are trivial. In $%
G,$ convergence to $id$ is determined by pointwise convergence to $id$ of
the coordinates in $\Pi _{G_{n}}K_{n}$ as shown in Theorem \ref{gcomp}.

\subsection{There exist wild braids in $\overline{B_{\infty }}$}

Elements of the finite braid groups $B_{n}$ can be understood as path
homotopy classes of homeomorphisms of the plane which start at the identity,
and place the finite set $F_{n}$ onto itself at the last moment of the
isotopy$.$ However $\overline{B_{\infty }}$ cannot be understood in
corresponding fashion due to the existence of `wild' braids.

By definition elements of the pure braid subgroup $G\subset \overline{%
B_{\infty }}$ can be understood as equivalence classes of paths of
embeddings of $Z^{+}$ into the plane which start and stop at $id_{Z^{+}}.$
The set $Z^{+}$ is allowed to move off itself in midjourney. Consider for $%
t\in \lbrack 0,1]$ the following isotopy $g_{t}$ of $Z^{+}.$ For $t\in
\lbrack 1-\frac{1}{n},1-\frac{1}{n+1}]$ let $(n,0)$ perform a small
clockwise orbit once around $(n+1,0)$ while leaving all other points of $%
Z^{+}$ fixed. Let $S$ denote a simple closed curve bounding a convex disk $D$
such that $\{(1,0),(2,0)\}\subset int(D)$ and $\emptyset =D\cap
\{(3,0),(4,0),...\}.$ For $t\in \lbrack 0,1)$ if we extend the isotopy $%
g_{t} $ to include $S$ we notice the diameter of $S$ is forcibly stretched
by arbitrarily large amounts. Hence the isotopy cannot be extended to $S$
for $t\in \lbrack 0,1]$ and hence the braid determined by $g_{t}$ is not
determined by an ambient isotopy of the plane.

\end{document}